\documentclass[12pt, a4paper]{article}
\usepackage[utf8]{inputenc}
\usepackage[english]{babel}
\usepackage[OT1]{fontenc}
\usepackage{amsmath}
\usepackage{amsfonts}
\usepackage{amssymb}
\usepackage{amsthm}
\usepackage{tikz}
\usetikzlibrary{positioning, calc}
\usepackage{graphicx}
\usepackage{tikz}
\usetikzlibrary{positioning, calc}
\usepackage[labelformat=empty]{caption}
\usepackage{float}
\usepackage{hyperref}
\usepackage{amstext} 
\usepackage{array}

\newtheorem{prop}{Proposition}
\newtheorem{coro}{Corollary}
\newtheorem{statement}{Statement}

\newtheorem{conj}{Conjecture}

\title{Searching for groups related to pseudo-composition algebras}
\author{Vsevolod A. Afanasev}
\date{}

\begin{document}

\maketitle

\begin{abstract}
We study the class of idempotent-generated pseudo-composition algebras, which is a subclass of the family of axial algebras. More specifically, we utilise the group-algebra correspondence, natural to the axial framework in order to study some automorphism subgroups of such pseudo-composition algebras.
\end{abstract}

\section{Introduction}
\let\thefootnote\relax\footnotetext{This work was supported by the Russian Science Foundation, project 23-41-10003, https://rscf.ru/en/project/23-41-10003/}

The theory of axial algebras, started in its current generality by Hall, Rehren and Shpectorov \cite{Axial}, contains an interesting class of examples: idempotent-generated pseudo-composition algebras, i.e. commutative algebras equipped with a non-zero symmetric bilinear form $\varphi$ such that
$$x^3=\varphi(x,x)x$$
for any element of an algebra.

Various properties of these objects, relayed in \cite{Osb}, allowed the authors of \cite{GorMamStar} to develop a more axial view of their properties. Historically, axial algebras were very closely connected to groups: for instance, the Matsuo algebras, which form a large subclass of the so-called axial algebras of Jordan type, are constructed from $3$-transposition groups. 
For more information on axial algebras, we refer interested readers to the survey \cite{ShpecMcSurv}.

In that regard, the newly established view is very useful as the axial framework provides easy methods for answering various group-theoretic inquiries about pseudo-composition algebras.
Hence, in this article we study the group-algebra correspondence for this class and provide multiple statements that may be used to study the Miyamoto group of any finitely-generated axial pseudo-composition algebra. 

More generally, we aim to use the parameter-based classification obtained in \cite{GorMamStar} in order to start finding answers to the following

\textbf{Question.} Which groups can be realised as the Miyamoto groups of axial pseudo-composition algebras?

Naturally, this question is incredibly grand and initially we wish only to obtain some interesting examples, in order to notice potential patterns and restrictions that are imposed on such groups. Another goal we have is to find some small-dimensional pseudo-composition representations of small non-abelian finite simple groups.

In the first section we shall study the algebras, generated by two idempotents, namely $a$ and $b$, and determine, which groups can be realised as their Miyamoto groups.

We will then make use of these results in order to move on to the algebras generated by three idempotents, for which we take a twofold approach: both trying to find some groups explicitly as well as attempting to determine whether some specific groups are Miyamoto groups for some pseudo-composition algebras.

For most of the paper, the ground field is assumed to be $\mathbb{R}$, i.e. a field of real numbers, however we will sometimes consider finite fields, preemptively informing the reader of such a change.
\paragraph*{Acknowledgements.}
The author would like to thank A. Mamontov, A. Staroletov and an anonymous referee for their insightful comments and incredibly helpful suggestions.

\section{The case of $2$-generated algebras}

Here we shall consider the pseudo-composition algebras, generated by two idempotents $a$ and $b$. One need not assume that these idempotents are axes, however it is known that any primitive idempotent in a pseudo-composition algebra automatically satisfies all the axial properties \cite{Osb}.  More generally, idempotent-generated pseudo composition algebras are actually axial, thus any such pseudo-composition algebra $A$ can be decomposed into a direct sum of eigenspaces $A_\lambda(a)$ with the following fusion law
\begin{center}
    $$\begin{array}{|c|c|c|c|}
    \hline
    * & 1 & -1 & \frac{1}{2}\\ \hline
    1 & 1 & -1 & \frac{1}{2}\\ \hline
    -1 & -1 & 1 & \frac{1}{2} \\ \hline
    \frac{1}{2} & \frac{1}{2} & \frac{1}{2} & -1,1 \\ \hline

    \end{array}
    $$

    Fusion law for $\mathcal{PC}(-1)$ -- a class of pseudo-composition algebras.
\end{center}

Additionally, this fusion law is naturally $\mathbb{Z}_2$-graded with $A_{\frac{1}{2}}$ being the "odd"\ part of the grading. Using these notions we can, given an axis $a$, construct an order $2$ automorphism $\tau_a$ of $A$: any element $b \in A$ can 
 be presented as follows:
$$b=b_+ + b_{-},$$
where $b_{-} \in A_{\frac{1}{2}}(a),$ while $b_{+} \in A_{-1}(a) \oplus A_{1}(a).$
Now denote
$$\tau_a(b)=b_+- b_{-}.$$
This automorphism is called the Miyamoto involution.
We then define the Miyamoto group of $A$ as the group generated by $\{\tau_a|\ a \text{ is an axis in A }\}$ . Note that generally this needs not be the whole automorphism group of $A$. 

Per this axis $\to$ automorphism correspondence, it is clear that the Miyamoto subgroups for the $2$-generated algebras $\langle\langle a, b \rangle\rangle$, are $D_n$ where $n=|\tau_{ab}|=|\tau_a\tau_b|$. Therefore to classify them we only need to find the possible orders of $\tau_{ab}$. By \cite{GorMamStar}, this matrix depends on the parameter $\alpha=(a, b)$, where $(\cdot, \cdot) $ is the Frobenius form, corresponding to the axial algebra (for pseudo-composition algebras it coincides with the bilinear form $\varphi$ from the definition).

The matrix may be written explicitly as follows (we assume that linear transformations act on the left)

$$\tau_{ab}=\frac{1}{9}\left( 
\begin{array}{ccc} 
64 \alpha^2-16\alpha-3 & 24\alpha & 32\alpha^2 +4\alpha-6 \\ 
8-8\alpha & -3 & -4\alpha -2 \\
-32\alpha -4 & -12 & -16\alpha +1\\
\end{array}  \right)$$
here the basis is $a,\ b,\ ab$.

\begin{prop}\label{2gend}
Let $l \in \mathbb{N}$. Then the following correspondence holds for the field of real numbers $\mathbb{R}$:

$$\begin{array} {|c|c|}
\hline
\text{Order of }\tau_{ab} & \alpha \\ \hline
2 & 1/4\\ 
3 & -1/8,\ 5/8\\ 
4 & 1/4 \pm 3/(4\sqrt{2})\\ 
5 & 7/16\pm 3\sqrt{5}/16,\ 1/16 \pm 3\sqrt{5}/16 \\
6 & 1/4\pm (3\sqrt{3})/8 \\
7 & \text{real root of } x^3-\frac{3}{8}x^2-\frac{9}{32}x+\frac{13}{512} \text{ or } x^3-\frac{9}{8}x^2+\frac{3}{32}x+\frac{43}{512}\\
8 & 1/8\left(2 \pm 3 \sqrt{2 \pm \sqrt{2}} \right) \\
9 & \text{real root of } x^3-\frac{3}{4}x^2-\frac{15}{64}x+\frac{19}{512} \text{ or }
x^3-\frac{3}{4}x^2-\frac{15}{64}x+\frac{73}{512} \\
10 & 1/16\left(4 \pm 3 \sqrt{2(5\pm \sqrt{5})}\right) \\
11 & \vtop{\hbox{\strut$\text{real root of } x^5-\frac{7}{8}x^4-\frac{5}{16}x^3+\frac{127}{512}x^2+\frac{119}{4096}x-\frac{263}{32768}$}\vspace{1mm} \hbox{\strut$ \text{ or } x^5 -\frac{13}{8}x^4+\frac{7}{16}x^3+\frac{145}{512}x^2-\frac{337}{4096}x-\frac{197}{32768}$}} \\
12 & 1/8 \left(2 \pm 3 \sqrt{2 \pm \sqrt{3}} \right) \\
\hline

\end{array}$$

\end{prop}

The proposition can be proved as follows: for a given order $k$ consider the matrix $$T=\tau_{ab}^k-I.$$ 
Then one can find the greatest common divisor of $T_{i,j}$ and check, whether it has any real roots. This procedure can easily performed via almost any modern computer algebra system and can be easily done for orders larger than $12$.

\textbf{Remark.} It can be seen that Proposition~\ref{2gend} holds for any finite extension of $\mathbb{Q}$, that contains the necessary square roots.

We have shown $|\tau_{ab}|$ for algebras over $\mathbb{R}$ can have many values and it can be conjectured that in this case it can actually take on any value. Another way to obtain the desired orders is to pass to the field of positive characteristic.
For instance, in the case of $\mathbb{F}_7$ we have the matrix looking as follows
$$\overline{\tau}_{ab}=\frac{1}{\Bar{2}}\left( \begin{array}{ccc}
 \alpha^2- \Bar{2}\alpha-\Bar{3}& \Bar{3}\alpha  & \Bar{4}\alpha^2+\Bar{4}\alpha - \Bar{6} \\
  \Bar{1}-\alpha   & -\Bar{3} & -\Bar{4}\alpha - \Bar{2} \\
  -\Bar{4} \alpha - \Bar{4}& -\Bar{5} & -\Bar{2}\alpha + \Bar{1}\\
     
\end{array} \right).$$

Substituting values of $\alpha$ from this field yields the following orders:
\begin{center}
$\begin{array}{|c|c|c|c|c|c|c|c|}
\alpha &\Bar{0} & \Bar{1} & \Bar{2} & \Bar{3} & \Bar{4} & \Bar{5} & \Bar{6}\\ \hline
|\overline{\tau}_{ab}|& 4 & 7 & 2 & 7 & 4 & 3 & 3 \\ 
\end{array}.$

\end{center}

Similarly, for $\mathbb{F}_{11}$ we have 
\begin{center}
$\begin{array}{|c|c|c|c|c|c|c|c|c|c|c|c|}
\alpha & \Bar{0} & \Bar{1} & \Bar{2} & \Bar{3} & \Bar{4} & \Bar{5} & \Bar{6} &
\Bar{7} &\Bar{8} &\Bar{9} &\Bar{10}\\ \hline
|\overline{\tau}_{ab}|& 5 & 11 & 3 & 2 & 3 & 11 & 5 & 5 & 6 & 6 &5\\ 
\end{array}.
$
\end{center}

It can thus be seen that many orders of $\tau_{ab}$ are reachable for the case of zero characteristic, however working over positive characteristic can also be useful in that regard.

\section{$3$-generated algebras}

In \cite{GorMamStar}, the authors construct $3$-generated axial pseudo-composition algebras $A=\langle \langle a,b,c \rangle \rangle$ of dimension $8$ with $4$ parameters $\alpha=(a,b),\ \beta=(b,c),\ \gamma=(a,c)$ and $\psi=(b,ac)$ defining each algebra's multiplication table and the Gram matrix. These algebras are distinguished by the following property: any other $3$-generated axial pseudo-composition algebra with the similar Frobenius form and parameter values is a homomorphic image of $A$.

We will call such algebras \textit{universal} axial pseudo-composition algebras.

The goal and strategy of this section are similar: we intend to construct the Miyamoto involutions, corresponding to the generators and manipulate the parameter values in an attempt to discover interesting Miyamoto groups.
Said matrices will now depend on four aforementioned parameters.

We will be constructing the matrices in the following basis
$$a,\ b,\ c,\ ab,\ bc,\ ac,\ a(bc),\ b(ac).$$

Given that $\tau_a(b)=1/3(8\alpha a- b-4ab)$, and $\tau_a(c)=1/3(8\gamma a-c-4ac)$ we construct the matrices via the automorphism multiplicative property:

\vspace{2mm}

$\tau_a=\left( 
\begin{array}{llllllll}
 1 & 8/3\alpha & 8/3\gamma & 4/3\alpha&  8/3\psi & 4/3\gamma &4/3 \psi &  8/3\alpha\gamma+\beta-2/3\psi\\
 0 & -1/3 & 0 &-2/3 & 0 & 0 & 0 & -1/3\gamma\\
  0 & 0 & -1/3 & 0 & 0 & -2/3 & 0 & \alpha\\
  0 & -4/3 & 0 & 1/3 & 0 & 0 & 0& 2/3\gamma\\ 
  0 & 0 & 0 & 0 & -1/3 & 0 & -2/3 & 1/3\\
   0 & 0 & -4/3 & 0 & 0 & 1/3 &0 & -2\alpha  \\
   0 & 0 & 0 & 0 & -4/3 & 0 & 1/3 & -2/3  \\
   0 & 0 & 0 & 0 & 0 & 0 & 0 & -1  \\
  
\end{array}
\right).$

\vspace{2mm}

In similar fashion we derive that 

$$\tau_b=\left(
\begin{array}{cccccccc}
     -1/3& 0& 0& -2/3& 0& 0& -1/3\beta& 0 \\
      8/3\alpha& 1& 8/3\beta& 4/3\alpha& 4/3\beta& 8/3\psi& 8/3\alpha\beta+\gamma-2/3\psi& 4/3\psi\\
      0& 0& -1/3& 0& -2/3& 0& \alpha& 0\\
      -4/3& 0& 0& 1/3& 0& 0& 2/3\beta& 0 \\
       0& 0& -4/3& 0& 1/3& 0& -2\alpha& 0\\
        0& 0& 0& 0& 0& -1/3& 1/3& -2/3\\
        0& 0& 0& 0& 0& 0& -1& 0 \\
0& 0& 0& 0& 0& -4/3& -2/3& 1/3 \\
\end{array}
\right), $$

$$\tau_c=\left(
\begin{array}{llllllll}
 -1/3& 0& 0& -4/3\beta& 0& -2/3& -\beta& 1/3\beta \\
 0& -1/3& 0& -4/3\gamma& -2/3& 0& 1/3\gamma& -\gamma\\
  8/3\gamma& 8/3\beta& 1& -4/3\alpha+8/3\psi& 4/3\beta& 4/3\gamma& \xi& \xi\\
  0& 0& 0& -1/3& 0& 0& 1/3& 1/3 \\
  0& -4/3& 0& 0& 1/3& 0& -2\gamma& 2/3\gamma \\
  -4/3& 0& 0& 0& 0& 1/3& 2/3\beta& -2\beta\\
  0& 0& 0& 4/3& 0& 0& -1/3& 2/3\\
  0& 0& 0& 4/3& 0& 0& 2/3& -1/3\\
  
\end{array}
\right),
$$
where $\xi=8/3\beta\gamma+1/3\alpha-2/3\psi$.

We will now try finding the optimal parameter values so that the products of obtained matrices are of finite order and investigate the groups generated by them. In this paper we will be mainly focused on obtaining finite groups, however infinite Miyamoto groups are interesting to study as well.

It is clear that the results found in the previous section may only be of partial service: while we know that, for instance $\tau_{ab}^k$ acts trivially on $\langle a,b \rangle$ for a certain value of $\alpha$, we do not know its behaviour on the rest of $A$. Therefore, the best we can hope for is that some multiple of the order found for the 2-generated case is the order of $\tau_{ab}$. 

\begin{prop}\label{oneparam}
    The characteristic polynomial of  $\tau_{ab},\ \tau_{bc},\ \tau_{ac}$ lies in  $\mathbb{R}(\alpha),\mathbb{R}(\beta)$ and $\mathbb{R}(\gamma)$ respectively.
\end{prop}

We shall first provide the proof for the universal $3$-generated axial pseudo-composition algebra $A$, dealing with the quotients later.

It suffices to consider $\tau_{ab}$, since other cases are reduced to this one by permuting the basis vectors and (for $\tau_{ac}$) replacing $b(ac)$ with $c(ab)$ in the basis (clearly not changing any invariants and thus preserving the order). 
Multiplying $\tau_a$ and $\tau_b$ we get

\vspace{2mm}

\small
\noindent
$
\displaystyle
\begin{pmatrix}
 & & & & R_1 & & &\\
  -\frac{8}{9}\alpha+\frac{8}{9} & -\frac{1}{3} & -\frac{8}{9}\beta & -\frac{4}{9}\alpha-\frac{2}{9} & -\frac{4}{9}\beta & \frac{4}{9}\gamma-\frac{8}{9}\psi & -\frac{8}{9}\alpha\beta-\frac{4}{9}\beta-\frac{1}{9}\gamma+\frac{2}{9}\psi & -\frac{1}{9}\gamma-\frac{4}{9}\psi \\ 
0 & 0 & \frac{1}{9} & 0 & \frac{2}{9} & -\frac{4}{3}\alpha+\frac{2}{9} & -\alpha-\frac{2}{9} & \frac{1}{3}\alpha+\frac{4}{9} \\ 
 -\frac{32}{9}\alpha-\frac{4}{9} & -\frac{4}{3} & -\frac{32}{9}\beta & -\frac{16}{9}\alpha+\frac{1}{9} & -\frac{16}{9}\beta & -\frac{8}{9}\gamma-\frac{32}{9}\psi & -\frac{32}{9}\alpha\beta+\frac{2}{9}\beta-\frac{16}{9}\gamma+\frac{8}{9}\psi & \frac{2}{9}\gamma-\frac{16}{9}\psi \\  0 & 0 & \frac{4}{9} & 0 & -\frac{1}{9} & -\frac{4}{9} & \frac{2}{3}\alpha+\frac{4}{9} & \frac{1}{9} \\ 
0 & 0 & \frac{4}{9} & 0 & \frac{8}{9} & \frac{8}{3}\alpha-\frac{1}{9} & \frac{1}{9} & -\frac{2}{3}\alpha-\frac{2}{9} \\ 
0 & 0 & \frac{16}{9} & 0 & -\frac{4}{9} & \frac{8}{9} & \frac{8}{3}\alpha+\frac{1}{9} & -\frac{2}{9} \\ 
0 & 0 & 0 & 0 & 0 & \frac{4}{3} & \frac{2}{3} & -\frac{1}{3} 

\end{pmatrix}
$

\vspace{2mm}

Where $$R_1=(\frac{64}{9}\alpha^2-\frac{16}{9}\alpha-\frac{1}{3},\ \frac{8}{3}\alpha,\  \frac{64}{9}\alpha\beta-\frac{8}{9}\gamma-\frac{32}{9}\psi,\  \frac{32}{9}\alpha^2+\frac{4}{9}\alpha-\frac{2}{3},\  \frac{32}{9}\alpha\beta-\frac{16}{9}\gamma+\frac{8}{9}\psi,$$ $$-\frac{32}{9}\alpha\gamma+\frac{64}{9}\alpha\psi-\frac{4}{3}\beta-\frac{4}{9}\gamma+\frac{8}{9}\psi,\frac{64}{9}\alpha^2\beta+\frac{8}{9}\alpha\beta+\frac{32}{9}\alpha\gamma-\frac{64}{9}\alpha\psi-\beta+\frac{4}{9}\gamma-\frac{8}{9}\psi,\  \frac{8}{9}\alpha\gamma+\frac{32}{9}\alpha\psi+\frac{1}{3}\beta-\frac{8}{9}\gamma-\frac{2}{9}\psi) $$
\normalsize
is the first row of the matrix, put in such form for compactness reasons.
The characteristic polynomial of this matrix turns out as follows:
$$x^8-\frac{8}{9}(8\alpha^2+2\alpha+1)x^7+\frac{4}{27}(256\alpha^3-48\alpha^2-24\alpha+5)x^6-\frac{4}{81}(1024\alpha^4+512\alpha^3-336\alpha^2-100\alpha+34)x^5+$$
$$+\frac{2}{81}(4096\alpha^4-1024\alpha^3-192\alpha^2+32\alpha-77)x^4-\frac{4}{81}(1024\alpha^4+512\alpha^3-336\alpha^2-100\alpha+34)x^3+$$
$$+\frac{4}{27}(256\alpha^3-48\alpha^2-24\alpha+5)x^2-\frac{8}{9}(8\alpha^2+2\alpha+1)x+1.$$

Any $3$-generated axial pseudo-composition algebra is a quotient of the universal $3$-generated axial pseudo-composition algebra by the latter's definiton. So we can consider the algebra $A_1=A/I$.

The ideal $I$ is generated as a vector space by the elements $v_1\ldots v_k$. We add vectors $e_1\ldots e_{p}$ to this basis, so that $k+p=8$ and so consider the basis $e_1,\ldots e_{p},\ v_1,\ldots v_k$. Since by \cite[Corollary 3.11]{KhasMcShpec}, we have that $I$ is invariant under the action of the Miyamoto group, and, in particular $\tau_{ab}$, we can present the matrix for $\tau_{ab}$ in this basis, as a following block matrix:

$$\tau_{ab}=\left(\begin{array}{cc}
    \left[\tau_{ab}\right] & * \\
     0 & I
\end{array}\right),$$
where $I$, by abuse of notation, denotes the elements $v_i$ under $\tau_{ab}$. It is clear that after taking the quotient only $\left[\tau_{ab}\right]$ will remain and it is known from linear algebra, that the characteristic polynomial of $\tau_{ab}$ is the multiple of characteristic polynomials of the diagonal blocks, thus the characteristic polynomial of the automorphism of $A_1$ corresponding to $\tau_{ab}$ divides $\chi_{\tau_{ab}}$ and thus also lies in $\mathbb{R}[\alpha]$. The proposition is proved.

As an easy consequence of this we have

\begin{coro}\label{divi}
    Assume $A_1=\langle \langle a_1,\  b_1,\ c_1 \rangle \rangle$ is a pseudo-composition algebra, generated by $3$ idempotents such that $\varphi(a)=a_1,\ \varphi(b)=b_1,\  \varphi(c)=c_1$, where $\varphi: A \to A_1$ -- a natural homomorphism by the ideal $I$. 
    Then $|\tau_{a_1b_1}|$ divides $|\tau_{ab}|$.

\end{coro}

We will provide another technical statement.

\begin{statement}
The elements of the first, second and fourth column of $\tau_{ab}^n$ are either $0$ for column index $j$ such that  $j \in \{3,\ 5,\ 6,\ 7, \ 8 \}=J$ or non-zero and lie in $\mathbb{R}[\alpha]$ otherwise.
\end{statement}

This can be proved by obvious induction on $n$. The base clearly holds, so we proceed to the induction step by considering $n=k+1$. By induction, the only non-zero elements of the considered columns are $a_{i,1}\ a_{i,2},\ a_{i,4}$, thus proving the zeroes claim.

The second claim follows from the fact that, while the first, second and fourth rows are the only ones containing parameters other than $\alpha$, by induction on $n$, they contain them only in places $a_{i,j}$ for $j \in J$, while $i=1,\ 2,\ 4$, so they cancel out upon multiplication. The statement is thus proved. $\square$

This effect can be easily explained by the fact that any $\tau_{ab}$ leaves any vector $v \in \langle \langle a,b \rangle \rangle$ in this subalgebra, i.e. $\tau_{ab}^n v \in \langle a, b, ab \rangle$.

Because of the previous statement, we are able to provide a simple yet useful and somewhat universal method for determining suitable parameters (or proving their nonexistence) for a given degree of $\tau_{ab}$. Consider the elements $f_1(\alpha)=(\tau_{ab}^k)_{1,1}-1$ and $f_2(\alpha)=(\tau_{ab}^k)_{1,2}$. These are polynomial in $\alpha$ and contain no entries of $\beta,\  \gamma$ or $\psi$ (this can be confirmed by looking at the matrix). It is clear that if there exists $\alpha_0$ such that $\tau_{ab}^k(\alpha_0)=I$ then $$f_1(\alpha_0)=f_2(\alpha_0)=0. \eqno{(**)}$$

Note that one can use any non-zero element from the first,\ second or fourth column of the matrix as additional equations or as replacements for $f_1$ or $f_2$ and then employ any suitable real root-finding algorithm for determining $\alpha_0$.

Returning to the Miyamoto groups we see that, for the $3$-generated case, the Miyamoto involutions $\tau_{a}$ and $\tau_{b}$ commute only in rather specific situations.

\begin{prop}\label{commutoperat}
Assume $a,\ b$ and $c$ are axes, generating a pseudo composition algebra $A^\prime$ and, without the loss of generality, that $|\tau_{ab}|=2$. Then $dim(A^\prime) \leq 4$.
\end{prop}

We start our proof by again working in the universal algebra $A$.

In this case we have
$$f_1=\frac{1}{81}(4096\alpha^4-3072\alpha^3-576\alpha^2+464\alpha+33)-1,\ f_2=\frac{1}{27}(512\alpha^3-256\alpha^2-64\alpha+24)$$

The only common factor of these polynomials is $(\alpha-\frac{1}{4})$ with $\frac{1}{4}$ setting $|\tau_{ab}|=4$.

Now we move on to the quotient algebras. Consider the ideal $I \subseteq A$, generated by elements $\tau_{ab}^2x_i-x_i$ where $x_i$ are all the basis vectors outlined at the begining of the section.

We provide these elements here in the form of a matrix (the entries are multiplied by $9$ for compactness reasons, so the matrix actually contains images of $\tau_{ab}^2 9x_i-9x_i$)

$$\left(
\begin{array}{cccccccc}
0& 0& 0& 0& 0& 0& 0& 0 \\ 
0& 0& 0& 0& 0& 0& 0& 0 \\ 
-\frac{16}{3}\beta+\frac{32}{3}\gamma-\frac{64}{3}\psi& \frac{32}{3}\beta-\frac{16}{3}\gamma-\frac{64}{3}\psi& -16& -\frac{64}{3}\beta-\frac{64}{3}\gamma+\frac{128}{3}\psi& 8& 8& 16& 16 \\ 
0& 0& 0& 0& 0& 0& 0& 0 \\ 
-\frac{20}{3}\beta-\frac{32}{3}\gamma+\frac{64}{3}\psi& \frac{40}{3}\beta-\frac{20}{3}\gamma-\frac{8}{3}\psi& 1& -\frac{8}{3}\beta+\frac{64}{3}\gamma-\frac{128}{3}\psi& -14& 4& 8& 8 \\ 
-\frac{20}{3}\beta+\frac{40}{3}\gamma-\frac{8}{3}\psi& -\frac{32}{3}\beta-\frac{20}{3}\gamma+\frac{64}{3}\psi& 1& \frac{64}{3}\beta-\frac{8}{3}\gamma-\frac{128}{3}\psi& 4& -14& 8& 8 \\
-\frac{16}{3}\beta-\frac{4}{3}\gamma+\frac{44}{3}\psi& -\frac{4}{3}\beta+\frac{38}{3}\gamma-\frac{64}{3}\psi& \frac{1}{2}& \frac{44}{3}\beta-\frac{64}{3}\gamma+\frac{56}{3}\psi& 2& 2& -14& 4 \\ 
\frac{38}{3}\beta-\frac{4}{3}\gamma-\frac{64}{3}\psi& -\frac{4}{3}\beta-\frac{16}{3}\gamma+\frac{44}{3}\psi& \frac{1}{2}& -\frac{64}{3}\beta+\frac{44}{3}\gamma+\frac{56}{3}\psi& 2& 2& 4& -14 \\
\end{array}\right),
  $$
where the rows correspond to the elements of the basis.

It is clear that any ideal, such that the image of the automorphism has order $2$ on the quotient algebra, must contain $I$. 
We then note that the four last columns can be reduced to the following row echelon form:

$$\left(
\begin{array}{cccc}
8&8&16&16 \\ 
0& 18 & 36 & 36 \\
0&0&36&36 \\
0&0&0&18\\
0& 0& 0& 0 \\ 
0& 0& 0& 0 \\ 
0& 0& 0& 0 \\ 
\end{array}
\right).
$$

This implies that the span of the vectors $\tau_{ab}^2x_i-x_i$ is $4$-dimensional, meaning that $dim\  I \geq 4 $.

 Setting $\beta=1,\ \gamma=1/4=\psi$ we get that $dim(I)=4$ and $\tau_{ab}$ has order $2$ on $A/I$; keeping the same parameters and setting $b=c$ yields a $5$-dimensional ideal, which contains $I$, while setting $ab=0=c$ we get a $6$-dimensional $I$.

\begin{statement}[c.f. proposition 1]\label{3gend}
For the axial pseudo-composition algebra $A$, finitely generated by at least three axes, matrix $\tau_{ab}$ has its order divide $3,4,5,6,10$ for values of $ -1/8, 1/4,\frac{1}{16}\pm \frac{3\sqrt{5}}{16}, 5/8, \frac{7}{16}\pm \frac{3\sqrt{5}}{16}$ for $\alpha$. The same values work for $\tau_{ac}$ and $\gamma$ and $\tau_{bc}$ and $\beta$.
\end{statement}

\begin{proof}
The result for algebras, generated by $3$ axes can be directly confirmed via working with the universal algebra and referring to Corollary \ref{divi}. For the finitely generated case, consider $\tau_{ab}$ for axes $a,\ b \in A$ such that $\alpha=(a,b)$ is equal to some value from the result's statement. 
Since it is known that any axial pseudo-composition algebra contains a basis consisting of primitive axes \cite[Proposition 2]{GorMamStar}, we can choose any axis from the basis that is not equal to $a$ or $b$ as $c$ and consider the subalgebra $\langle \langle a, b, c \rangle \rangle$. By the definition this is a homomorphic image of the universal axial pseudo-composition algebra $\langle \langle a^\prime, b^\prime, c^\prime \rangle \rangle$ with the same parameter values. 
We have that $\tau_{a^\prime b^\prime}^k$ (where $k$ is one of the orders from the result's statement) fixes $a^\prime,\ b^\prime,\ c^\prime$, 
thus being of order $k$ for the universal algebra, and again by Corollary \ref{divi}, we have that $\tau_{ab}^{k_0}$, such that $k_0$ divides $k$, fixes $a,\ b,\ c$. 

By this argument $k_0$ depends on the initial choice of an axis $c$. Note however, that since $\alpha$ is fixed, $k$ does not depend on that choice so we can pick $k_0$ as the biggest divisor of $k$ that appears this way for some axis $c$, hence $\tau_{ab}^{k_0}$ will fix the whole algebra.
\end{proof}

This statement is another demonstration of the following effect: The parameter of a $2$-generated pseudo-composition algebra $\langle \langle a_1, a_2 \rangle \rangle$ with $a_i$ being axes, defines the order of the corresponding automorphism $\tau_{a_1a_2}$, but only up to a scalar multiple.

We will now focus our attention on finding groups, that can be generated by three automorphisms.

For example, by setting $\alpha=\beta=\gamma=-1/8$ we obtain a $3$ generated group, that is known to be a group of motions of the euclidean plane, generated by the reflections in the sides of an equilateral triangle. This group is known to be infinite, however it is possible to pass to its finite quotient by adding a relator $(\tau_a\tau_b^{\tau_c})^m$. More specifically:

\begin{prop}
Let $G=\langle x,y,z \rangle$ be the group generated by involutions $x,\ y,\ z$, such that $|xy|=|yz|=|xz|=3$. Additionally, let $k=|xy^z|$. Then $G$ is a quotient of $k^2 \rtimes D_6$.
\end{prop}
Here $k^2$ denotes the direct product of two copies of a cyclic group of order $k$, while $D_6$ is the dihedral group of order $6$.

We will provide the proof of this proposition for the sake of completeness:

\textit{Proof.}
Note that
$$[xy^z,yz^x]=y^zxxzxyxy^zyz^x=zyxyxy^zyz^x=zxyzyzyz^x=x^zz^x=(zx)^3=e.$$

Therefore the subgroup $K=\langle xy^z, yz^x \rangle$ is abelian (being isomorphic to $k^2$). It is also easy to see that it is a normal subgroup of $G$. Now consider the quotient $G/K$. In it, $\Bar{y}=\overline{z^x}$, therefore the quotient is generated by 2 involutions, making it a Dihedral group of order $2|xy|$, which is equal to $6$. Since the groups $K$ and $\langle x, z \rangle$ intersect trivially, we conclude that $G$ is a semidirect product of $K$ and $D_6$, concluding the proof.
$\square$

As a result we arrive at another family of groups, that may be realised as Miyamoto subgroups of pseudo-composition algebras. Now we only need to find the exact possible orders $\tau_a\tau_b^{\tau_c}$ can take. 
This can be accomplished by choosing the correct values of $\psi$, so the situation is somewhat similar to the previous cases, since we only have a single parameter.

Now, the minimal polynomial of this matrix, computed in GAP \cite{GAP}, looks as follows:
$$x^5+\frac{1}{81}(-1024\psi^2+32\psi+101)x^4+\frac{1}{243}\left(-\frac{32768}{3}\psi^3-1024\psi^2+1120\psi+\frac{404}{3}\right)x^3+$$
$$+\frac{1}{243}\left(-\frac{32768}{3}\psi^3-1024\psi^2+1120\psi+\frac{404}{3}\right)x^2+\frac{1}{81}(-1024\psi^2+32\psi+101)x-1$$
with one of its roots being equal to 1 and the other roots being

$$r_{1,2}=\frac{1}{9}\left(-(16\psi+2)\pm \sqrt{256\psi^2+64\psi-77}\right),$$

$$r_{3,4}=\frac{1}{162}\left(1024 \psi^2+256\psi-146 \pm \sqrt{-(1024\psi^2+256\psi-146)} \right).$$

One can then manually use these in order to see the following 

\begin{statement}
For $\alpha=\beta=\gamma=-1/8$ the values $\frac{5}{32}, -\frac{1}{8}, \frac{1}{64}-\frac{(9\sqrt{5})}{64}, -\frac{13}{32}$ of $\psi$ set the order of $ab^c$ to a divisor of $3,4,5,6$ respectively. 
\end{statement}

Note that after fixing the orders of involution products, $\psi$ is the sole remaining parameter and as such it may be insufficient to obtain a finite group. We will illustrate that by the following example:

Consider the following presentation of $A_5$
$$\langle x,y,z| x^2,y^2, z^2,(xy)^3,(yz)^3,(xyz)^3 \rangle.$$
We will attempt to find suitable parameters $\alpha,\ \beta,\ \gamma$ and $\psi$ to obtain  $\langle \tau_a,\ \tau_b,\ \tau_c \rangle \simeq A_5$. Now set $x=\tau_a,\ y=\tau_b,\ z=\tau_c$. This automatically fixes $\alpha=\beta=-1/8$ (due to Proposition \ref{2gend} and Statement \ref{3gend} this is the only possible value).

Additionally, in this presentation the elements $xz$ and $xy^z$ must have order $5$. To set the former we use one of the values of $\gamma$ from Proposition \ref{2gend} that sets the order to $5$ in the universal axial pseudo-composition algebra, with the uniqueness of said parameters also being guaranteed by Proposition \ref{2gend} and Statement \ref{3gend}, while for $xy^z$ we have to find suitable values of $\psi$. Therefore if no finite groups are produced for parameters, guaranteeing that order, this group cannot be obtained as a Miyamoto group of an $8$-dimensional pseudo-composition algebra, generated by $3$ axes. 

We will use the following construction: consider the axis $b^{\tau_c}$, which we will denote $d$. It can be seen, that $d$ has $\tau_b^{\tau_c}$ as its Miyamoto involution and that $d$ is an image of an idempotent under the automorphism and thus $d$ is an idempotent itself. We then consider the subalgebra $\langle \langle a, d \rangle \rangle$, which has $\tau_{ad}=\tau_a\tau_b^{\tau_c}$. Therefore we need to only derive the values of $\psi$ from the equality $(a,d)=r$ with $r$ being one of the parameter values from Proposition \ref{2gend}, corresponding to $5$. If none of these parameters yield a finite group, our claim follows.

With these parameters (we pick $\gamma=\frac{1}{16}-\frac{3\sqrt{5}}{16}$, noting that the reasoning is exactly the same for its conjugate).

With this setup, the angle $(a,d)=-\frac{4}{3}\psi+\frac{3\sqrt{5}+1}{48}$. 
We need to find four values of $\psi$ such that $(a,d) \in \{ \frac{1}{16}\pm \frac{3\sqrt{5}}{16},\ \frac{7}{16}\pm \frac{3\sqrt{5}}{16}\}$.

Thus we get that $$\psi \in \{ -\frac{1+3\sqrt{5}}{32},\ \frac{-1+6\sqrt{5}}{32},\ -\frac{10+3\sqrt{5}}{32},\ \frac{-5+3\sqrt{5}}{16} \},$$
where the values correspond to $\frac{1}{16}+ \frac{3\sqrt{5}}{16},\ \frac{1}{16}- \frac{3\sqrt{5}}{16},\ \frac{7}{16}+ \frac{3\sqrt{5}}{16},\ \frac{7}{16}- \frac{3\sqrt{5}}{16}$ respectively.

Two main ways will be used to eliminate values of $\psi$. Either $|\tau_a\tau_d|=\left|\tau_a\tau_b^{\tau_c}\right|$ is indeed equal to $5$ and we use GAP to check, whether the group generated by Miyamoto involutions is finite or if the order is equal to $5k$ (usually $k=2$), in which case we check, whether the corresponding algebra is simple by computing the rank of its Gram matrix, as the algebra being simple implies there is no way to find a non-trivial quotient with the required order. 

The first and second values of $\psi$ are setting the order of $\tau_a\tau_d$ to $5$ but the group generated by involutions is infinite, while the latter two values set the order to $10$, with the corresponding algebras being simple. Hence our claim follows.

We have considered the presentation of $A_5$ that has $(xy)^3$, $(yz)^3$ and $(xz)^5$ as relators of product orders. Since all element orders for this group lie in the set $\{1,\ 2,\ 3,\ 5\}$, by Proposition \ref{commutoperat}, we only need to consider orders $3$ and $5$ as possible product orders. However all other presentations can be dealt with in a similar way as for the case $(xy)^3,\  (yz)^5,\ (xz)^5$, one of two elements $xz^y$ or $yz^x$ always has order $2$ or $3$ either eliminating this case due to Proposition \ref{commutoperat} or reducing them to the above case by replacing $z$ with its conjugate. Likewise it can be manually checked that presentations with all generator products having order $5$ set the order of $xy^z$ to $3$ and can be rewritten into the previously mentioned presentations.

The technique used here, while simple, is generalisable, allowing one to find potential parameters for the pseudo-composition algebra for a given group, provided said group's presentation consists of relators of the form $a_ia_j^w$ where $a_i,\ a_j$ are among the generating involutions and $w$ is a word from the group. In particular, this highlights another connection to the class of $n$-transposition groups, which have this exact type of relators by definition. 

Nonetheless, we cannot confirm the non-existence of a pseudo-composition algebra generated by $3$ axes for $A_5$ as there may exist quotient algebras with it as a Miyamoto subgroup. It is to be hoped that finer analysis of the Miyamoto group of a $4$-dimensional pseudo composition algebra generated by $3$ axes will shed light on the validity of this non-existence hypothesis.

However for the second smallest non-abelian simple group, $PSL(2,7)$, the situation is more positive:

It is known that $PSL(2,7)$ satisfies the presentation
$$\langle t,s| t^2=s^3=(st)^7=[s,t]^4\rangle. $$

This can be transformed into a presentation with three generating involutions $a=t,\ b=t^s$ and $c=t^{s^2}$. It can be derived that $|a c|=|ab|=|bc|=4$ and that $|ab^c|=3$. This allows us, by setting $\alpha=\beta=\gamma=\frac{1}{4}$ and $\psi=\frac{5}{32}$ to obtain $\langle \tau_{a},\tau_{b},\tau_{c}\rangle=PSL(2,7)$.
The rank of the Gram matrix for the corresponding algebra is equal to 8, meaning that the algebra is simple. 

Therefore our work with two smallest non-abelian finite simple groups can be summarised in the following

\begin{prop}
Let $A$ be the universal axial pseudo composition algebra generated by $3$ axes $a,\ b,\ c$ over the field $\mathbb{R}$. Then we following statements hold.
\begin{enumerate}
    \item The alternating group $A_5$ cannot be generated by Miyamoto involutions $\tau_a,\  \tau_b,\ \tau_c$;
    \item The group $PSL(2,7)$ can be generated by $\tau_a,\ \tau_b,\ \tau_c$ with parameter values $\alpha=\beta=\gamma=\frac{1}{4}$ and $\psi=\frac{5}{32}$.
\end{enumerate}

\end{prop}

\subsection*{Case of positive charactersitic}

Akin to the 2-generated case, we can try tweaking the characteristic of the ground field.

\begin{prop}\label{char5list}
    Assume $A=\langle \langle a,b,c \rangle \rangle$ is the universal axial pseudo-composition algebra defined over $\mathbb{F}_5$ with $a,\ b,\ c$ being  idempotents. Then the group generated by $\tau_a, \tau_b, \tau_c$ is either solvable or is isomorphic to $5^5:A_5,\ 5^5:~S_5,$  
    $5^2:~(5^2:~(SL(2,5):~2)),\ PSL(2,7),$ $PSL(3,5),  PSU(3,5), \ A_6,\ A_7$. 

    For the listed groups the Gram matrix of $A$ has rank $5$ for the first two cases and $8$ for all the other ones. 
\end{prop}

We note that while said groups appear for multiple combinations of parameters (with them not just being permutations of each other), the rank of the Gram matrix remains constant for all options.

We will provide a single combination of parameters for every group listed in the proposition. We will not be using the bar notation, even though the parameter values lie in $\mathbb{F}_5$:

$$
\begin{array}{|c|c|c|c|c|c|}
\hline
\alpha & \beta & \gamma & \psi & T=\langle \tau_a,\ \tau_b,\ \tau_c \rangle & rank(gram(A))\\ \hline
3 & 4 & 1 & 1 & 5^2:~(5^2:~(SL(2,5):~2)) & 4 \\
3 & 1 & 3 & 2 & 5^5:A_5 & 5\\
3 & 0 & 4 & 4 & 5^5:S_5 & 5\\
3 & 3 & 4 & 4 & PSL(2,7) & 8 \\
3 & 3 & 1 & 0 & PSL(3,5) & 8 \\
1 & 3 & 1 & 0 & PSU(3,5) & 8 \\
3 & 3 & 1 & 1 & A_6 & 8\\
3 & 3 & 1 & 4 & A_7 & 8\\

\hline

\end{array}
$$

The proposition can be confirmed by directly substituting various parameters into $\tau_a,\ \tau_b$ and $\tau_c$. 

We note that for solvable groups the rank of the Gram matrix takes all values from $1$ to $8$, excluding $6$. This suggests that said solvable groups may prove to be interesting examples if one wishes to attempt classifying the ideals inside such pseudo-composition algebras.

It is important to note that $T=\langle \tau_a,\ \tau_b,\ \tau_c \rangle$ is only a subgroup of the Miyamoto group.

\textbf{Example.}
Consider the case when $T\simeq A_6$. After using the parameter values from Proposition \ref{char5list} to the corresponding algebra $A$, one can find that $b(ac)$ is an idempotent in that algebra. By \cite[Proposition 2.7]{Osb}, it is also an axis, so we will denote it by $d$, thus having a corresponding Miyamoto involution, denoted by $\tau_d$. 

Applying the automorphism construction rules, which can be extracted by generalising \cite[Lemma 2.1]{GorMamStar}, we can construct $\tau_d$, which, per our definitions is also a generator in the Miyamoto group. However, as it turns out
$\langle \tau_a,\ \tau_b,\ \tau_c,\ \tau_d \rangle \simeq SL(8,5)$. 

In other words, the group $T$ can heavliy differ from the whole Miyamoto group of an algebra, which however should not discourage interest of studying such subgroups. We also note that the number of involutions in $SL(8,5)$ is significantly greater than the number of idempotents in $A$, thus showing that the idempotent $\leftrightarrow$ involution correspondence is not surjective.

Finally, we state the results of our reconnaissance into bigger values of positive characteristic.

After experimenting with some parameters over fields of bigger characteristic we have arrived at a

\begin{statement}
    The group $T=\langle \tau_a,\ \tau_b,\ \tau_c \rangle$ is isomorphic to $PSL(3,q)$ for the following combinations of $(q, \alpha, \beta, \gamma ,\psi)$:
    $$(7,\ \Bar{6},\ \Bar{2},\ \Bar{1},\ \Bar{5}),\; (11,\ \Bar{1},\ \Bar{1},\ \Bar{3},\ \Bar{0}),\; (13,\  \Bar{8},\ \Bar{10},\ \Bar{4},\ \Bar{2}).$$

    Similarly, $T$ is isomorphic to $PSU(3,q)$ for the following combinations of $(q, \alpha, \beta, \gamma ,\psi)$:
    $$(7,\ \Bar{6},\ \Bar{2},\ \Bar{1},\ \Bar{1}),\; (11,\ \Bar{1},\ \Bar{1},\ \Bar{3},\ \Bar{1}),\;  (13,\  \Bar{8},\ \Bar{10},\ \Bar{4},\ \Bar{1}).$$
\end{statement}

This prompts to present the following:

\begin{conj}
There exist parameters $\alpha,\ \beta,\ \gamma,\  \psi \in F$ such that $T=\langle \tau_a,\ \tau_b,\ \tau_c \rangle$ is isomorphic to $PSL(3,F)$ and a different set of parameters such that $T$ is isomorphic to $PSU(3,F)$, where $F$ is any prime field.
\end{conj}

This conjecture may be possible to tackle via constructing reasonably good presentations for these groups, akin to $PSL(2,7)$.

Concerning other questions, computing the full Miyamoto groups and automorphism groups for algebras with parameter sets we have already considered is certainly interesting, as it would allow us to deepen our understanding of said algebras and of automorphisms of axial algebras as a whole.

\end{document}